\documentclass[a4paper,12pt,reqno]{amsart}

\usepackage{amssymb,amsmath,amsfonts,amsthm,mathtools,dsfont}
\usepackage{latexsym,mathrsfs,pb-diagram}
\usepackage[pdftex]{graphicx}
\usepackage{pgf,tikz,tikz-cd}

\usepackage[hmargin=3.0cm,vmargin=2.5cm]{geometry}
\usepackage[plainpages=false,colorlinks,pdfpagelabels]{hyperref}
\hypersetup{
	urlcolor=black,
	citecolor=green,
	linkcolor=blue
}

\numberwithin{equation}{section}

\theoremstyle{definition}
\newtheorem{Def}{Definition}[section]

\theoremstyle{remark}

\newtheorem{Rem}[Def]{Remark}

\theoremstyle{plain}
\newtheorem{Prop}[Def]{Proposition}
\newtheorem{Cor}[Def]{Corollary}
\newtheorem{Thm}[Def]{Theorem}
\newtheorem{Lem}[Def]{Lemma}

\newcommand{\Z}{\mathbb{Z}}
\newcommand{\N}{\mathbb{N}}
\newcommand{\R}{\mathbb{R}}

\newcommand{\C}{\mathbb{C}}

\newcommand{\transp}[1]{\prescript{\mathrm{t} \!}{}{{#1}}}

\DeclareMathOperator{\ran}{\mathrm{ran}}

\newcommand{\dd}{\mathrm{d}}

\newcommand{\TT}{\mathbb{T}}

\newcommand{\D}{\mathscr{D}}

\newcommand{\cinfty}{\mathscr{C}^\infty}

\title{Global solvability for a class of pseudodifferential operators on the torus}
\author{I. A.~Ferra}
\address{Universidade Federal de São Carlos, Brazil}
\email{\texttt{igorferra@ufscar.br}}

\keywords{global solvability, pseudodifferential operators, periodic solutions.} 
\subjclass[2020]{35A01, 35B10, 58J40}
\begin{document}
\maketitle

\begin{abstract}
We give a complete characterization for the global solvability of a pseudodifferential operator
\[
P=D_t + c\left(t,D_x\right)
\]
on the $(N+1)$-dimensional torus $\TT^{N+1} = S^1_t\times \TT^N_x$. Our characterization is given in terms of diophantine conditions and a notion of super-logarithmic oscilation of the symbol of the imaginary part of $c\left(t,D_x\right)$.
\end{abstract}
\section{Introduction}

The local solvability of a classical pseudodifferential operator of principal type is completely characterized by the $(\Psi)$ condition of Nirenberg-Treves (see~\cite{dencker03}). In contrast, it was proved in~\cite{hounie79} that this condition is not necessary for the {\it global} solvability of vector fields on the $2$-dimensional torus $\TT^2$. Since then, some works have characterized the global solvability of classes of differential operators of the form $D_t + c\left(t,D_x\right)$ in the $(N+1)$-dimensional torus $S^1_t\times \TT^N_x$ (for instance~\cite{bcg21,bdg17}). In this work we present a complete characterization for the global solvability of
\begin{align}\label{eq:operator-P}
P = D_t + c\left(t,D_x\right),
\end{align}
where $c\left(t,D_x\right)$ is a pseudodifferential operators of order $m\in \R$ in the class introduced first in~\cite{cc17} for the real analytic case and in~\cite{fp23} for the smooth case that we recall here: a linear and continuous operator $Q\left(y,D_y\right):\cinfty\left(\TT^{n}\right)\to \cinfty\left(\TT^{n}\right)$ in the $n$-dimensional torus $\TT^n$ belongs to the class $S^m\left(\TT^{n}\times \Z^{n}\right)$ if given $\alpha\in\Z_+^n$, there exists $C_\alpha>0$ such that
\[
\left|D^\alpha_y Q\left(y,\eta\right)\right| \le C_\alpha\left(1+\left|\eta\right|\right)^{m}, \quad \forall y\in\TT^n, \eta\in\Z^n,
\]
where $Q\left(y,\eta\right)\doteq e^{-iy\eta}Q\left(y,D_y\right)e^{iy\eta}$ is the discrete symbol of $Q\left(y,D_y\right)$. We are assuming that the discrete symbol of our operator $c$ in $S^1_t\times \TT_x^N$ depends only on $(t,\xi) \in S^1\times \Z^N$ -- the so called {\it tube type} operator -- and it is worth to point out that we make no further assumptions about the class under study like being a sum of homogeneous operators. We write 
\[
c\left(t,D_x\right)=a\left(t,D_x\right)+ib\left(t,D_x\right)
\]
where $a\left(t,D_x\right)$ and $b\left(t,D_x\right)$ are the real and imaginary parts, respectively, of $c\left(t,D_x\right)$.

Our main result reveals two types of obstructions for the global solvability of $P$ given in~\eqref{eq:operator-P}:
\begin{enumerate}
	\item A diophantine condition on the averages of the symbol $c\left(t,\xi\right)$ -- see~\eqref{eq:DC};
	\item A super-logarithmic oscilation of primitives of $b\left(t,\xi\right)$ in intervals of size bigger than $\left|\xi\right|^{-m}$ -- see Definition~\ref{Def:condicao-alpha-beta}.
\end{enumerate}
The first obstruction is classic on the existence and regularity of global solutions of linear partial differential equations and it is the only point where the real part $a\left(t,D_x\right)$ contributes to the global solvability of $P$. The logarithmic growth of the symbol was explored before (see~\cite{agkm19,agk18}) for the problem of global hypoellipticity of a class of pseudodifferential operators on torus specially when it is possible to write $c\left(t,D_x\right)=c_0\left(t\right)Q\left(D_x\right)$.

It is interesting to note that our characterization is given only in terms of inequalities. Therefore, the conectedness of the certain sublevels of primitives of $b\left(t,\xi\right)$ -- a common property in previous works about this subject -- does not play an explicit role in our main result. When $c\left(t,\xi\right)$ is positively homogeneous in $\xi$, Corollary~\ref{Cor:homog-charact} shows how the inequalities recover previous results given in terms of the conectedness of sublevels.

To state our result we denote, for each $\xi\in\Z^N$, the averages
\begin{align*}
c_0(\xi) &= (2\pi)^{-1} \int_0^{2\pi} c(t,\xi)\dd t = a_0(\xi)+ib_0(\xi)
\end{align*}
%
%
%
and we consider the set
\[
\mathcal Z = \left\{\xi \in \Z^N : c_0(\xi) \in \Z\right\}.
\]
Given $\xi\in\Z^N$, $C\left(t,\xi\right)$ and $B\left(t,\xi\right)$ will denote primitives of $c(t,\xi)$ and $b\left(t,\xi\right)$, respectively, that sometimes will be chosen carefully depending on our pourposes (all the statements are independent of this choice). If $\xi \in \mathcal Z$ then $b_0\left(\xi\right)=0$, so $B(t,\xi)$ is a $2\pi$-periodic function on the first variable and we can view this function as a smooth function in $S^1$. In this case, $t_\xi \in S^1$ will denote a point where $B\left(t,\xi\right)$ assumes its maximum. Given two distinct points $t_1,t_2\in S^1$, we denote by $\left[t_1,t_2\right]\subset S^1$ the arc from $t_1$ to $t_2$ in the anticlockwise orientation. 
\begin{Def}
	\label{Def:condicao-alpha-beta}
	Let $D>0$.
	\begin{enumerate}
		\item We say that $P$ satisfies $(\alpha)_D^+$ at $\xi\in\Z^N$ if for every $t \in [0,2\pi]$ and every open interval $I^+\subset\left[t,t+2\pi\right]$ such that $\left|\xi\right|^{-m}\le \left|I^+\right|$, there exists $s \in I^+$ such that
		\[
		B(t,\xi)-B(s,\xi) \le D\log\left|\xi\right|.
		\]
		We say that $P$ satisfies $(\alpha)_D^-$ at $\xi$ if for every $t \in [0,2\pi]$ and every open interval $I^-\subset\left[t-2\pi,t\right]$ such that $\left|\xi\right|^{-m}\le \left|I^-\right|$, there exists $s \in I^-$ such that
		\[
		B(t,\xi)-B(s,\xi) \le D\log\left|\xi\right|.
		\]
		
		\item We say that $P$ satisfies the condition $(\beta)_{D}$ at $\xi \in \mathcal Z\setminus\{0\}$ when given $t \in S^1\setminus\{t_\xi\}$, one of the following conditions hold true:
		\begin{enumerate}
			\item  for every open interval $I^+\subset\left[t,t_\xi\right]$ in $S^1$ such that $\left|I^+\right|\ge \left|\xi\right|^{-m}$, there exists $s \in I^+$ such that
			\[
			B\left(t,\xi\right) - B\left(s,\xi\right) \le D\log\left|\xi\right|.
			\]
			\item for every open interval $I^-\subset\left[t_\xi,t\right]$ in $S^1$ such that $\left|I^-\right|\ge \left|\xi\right|^{-m}$, there exists $s \in I^-$ such that
			\[
			B\left(t,\xi\right) - B\left(s,\xi\right) \le D\log\left|\xi\right|.
			\]
		\end{enumerate} 
	\end{enumerate}
\end{Def}

The conditions $(\alpha)_D^\pm$ and $(\beta)_D$ do not depend on the choice of the primitive $B(t,\xi)$ of $b(t,\xi)$. Our main result is the following:
\begin{Thm}
	\label{Thm:main}
	The operator $P$ given in~\eqref{eq:operator-P} is globally solvable if and only if the following conditions hold true:
	\begin{enumerate}
		\item\label{Thm:main-item1}There are $C>0$ and $M\in\Z_+$ such that
		\begin{align}
			\label{eq:DC}
		\left|\tau+c_0(\xi)\right| \ge C(1+|\xi|)^{-M}, \quad \forall \tau \in \Z, \xi \in \Z^N\setminus\mathcal Z.
		\end{align}
		\item There exists $D>0$ such that:
		\begin{enumerate}
			\item\label{Thm:main-item2} Given $\xi\in\Z^N\setminus\mathcal Z$ such that $|\xi|\ge D$, $P$ satisfies $(\alpha)_D^+$ or $(\alpha)_D^-$ at $\xi$.
			\item\label{Thm:main-item3} Given $\xi\in\mathcal Z$ such that $|\xi|\ge D$, $P$ satisfies $(\beta)_D$ at $\xi$.
		\end{enumerate}
	\end{enumerate}
\end{Thm}
The work is organized as follows: in Section~\ref{Sec:pre} we set the basic results for the proof of Theorem~\ref{Thm:main}. In Section~\ref{Sec:nec} we present the necessity of the conditions of Theorem~\ref{Thm:main}, while Section~\ref{Sec:suf} is devoted to the sufficiency of such conditions. We finally show, in Section~\ref{Sec:app}, necessary conditions to the global solvability when the operators are positively homogeneous. In particular, Corollary~\ref{Cor:homog-charact} gives a complete characterization of the global solvability when $c\left(t,D_x\right)$ is positively homogeneous.
\section{Preliminaries}
\label{Sec:pre}

Consider, in $S^1$, the differential equation
\begin{align}\label{eq:ode}
\left(D_t + g(t)\right) u(t) = f(t),
\end{align}
where $f,g \in \cinfty\left(S^1\right)$. If we set $g_0=(2\pi)^{-1}\int_0^{2\pi}g(t)\dd t$ and if $G$ denotes a primitive of $g$ then:
\begin{enumerate}
	\item when $g_0 \in\C\setminus\Z$,~\eqref{eq:ode} has a unique solution given by
	\begin{align}\label{eq:solution-ode-outside-Z}
	u(t) &= \frac{i}{1-e^{-2\pi i g_0}}\int_0^{2\pi} e^{i\left(G(t-s)-G(t)\right)} f(t-s)\dd s\\ \notag
	&=\frac{i}{e^{2\pi i g_0}-1}\int_0^{2\pi} e^{i\left(G(t+s)-G(t)\right)} f(t+s)\dd s;
	\end{align}
	\item when $g_0 \in \Z$ and
	\[
	\int_{S^1} e^{iG(s)}f(s)\dd s = 0,
	\]
	the general solution of~\eqref{eq:ode} is
	\begin{align}\label{eq:solution-ode-inside-Z}
		u(t) = Ce^{-iG(t)}+i\int_{\left[t_0,t\right]} e^{i\left(G(s)-G(t)\right)} f(s)\dd s,
	\end{align}
	with $C \in \C$, where $t_0 \in S^1$ is arbitrary and $\left[t_0,t\right]$ denotes the arc from $t_0$ to $t$ in the anticlockwise orientation.
\end{enumerate}

It is important to note that we can choose any primitive $G$ of $g$ in~\eqref{eq:solution-ode-outside-Z} and~\eqref{eq:solution-ode-inside-Z}. Moreover, $e^{iG(t)} \in \cinfty\left(S^1\right)$ whenever $g_0 \in \Z$.
\begin{Def}
	We say that $P:\cinfty\left(\TT^{N+1}\right)\to \cinfty\left(\TT^{N+1}\right)$ is globally solvable when its range in closed in $\cinfty\left(\TT^{N+1}\right)$.
\end{Def}
By functional analysis arguments, $\overline{\ran P} = \left(\ker \transp P\right)^\bot$, that is $f \in \overline{\ran P}$ if and only if $\left\langle w, f \right\rangle =0$ for every $w \in \D'\left(\TT^{N+1}\right)$ such that $\transp P w = 0$, where $\transp P$ denotes the transpose of $P$. Our next result characterizes the set $\overline{\ran P}$ and its proof uses the fact that the transpose of $c\left(t,D_x\right)$ is given by $c\left(t,-D_x\right)$. This can be obtained by showing that
\[
\left\langle \transp c(t,D_x)\left(\varphi (t)\otimes \psi(x)\right), f(t)\otimes g(x)\right\rangle = \left\langle c(t,-D_x)\left(\varphi (t)\otimes \psi(x)\right), f(t)\otimes g(x)\right\rangle
\]
for every $f,\varphi  \in\cinfty\left(S^1\right)$ and $g,\psi \in \cinfty\left(\TT^N\right)$.
\begin{Prop}
	\label{Prop:ker-transpose}
	If $\xi\in\mathcal Z$ then $w_\xi\doteq e^{iC(t,\xi)}\otimes e^{-ix\xi} \in \cinfty\left(\TT^{N+1}\right)$ belongs to $\ker\transp P$, where $C\left(t,\xi\right)$ denotes a primitive of $c\left(t,\xi\right)$. Consequently, a function $f \in \cinfty\left(\TT^{N+1}\right)$ belongs to $\overline{\ran P}$ if and only if
	\begin{align}\label{eq:charac-kernel-transp}
		\int_{S^1} e^{iC(t,\xi)}\hat f(t,\xi)\dd t = 0, \quad \forall \xi \in \mathcal Z.
	\end{align}
\end{Prop}
\begin{proof}
Since $\transp P=-D_t + c\left(t,-D_x\right)$, it follows that
\begin{align*}
	\transp P w_\xi &= \left(-D_t + c(t,-D_x)\right) e^{iC(t,\xi)}\otimes e^{-ix\xi}\\
	&= \left(-c(t,\xi)+c(t,\xi)\right) e^{iC(t,\xi)}\otimes e^{-ix\xi}\\
	&= 0.
\end{align*}
If $f \in \overline{\ran P}$ and $\xi \in \mathcal Z$, then
\begin{align*}
	0&= \left\langle w_\xi,f\right\rangle = \left(2\pi\right)^{N}\int_{S^1} e^{iC(t,\xi)}\hat f(t,\xi)\dd t.
\end{align*}
Conversely, suppose that~\eqref{eq:charac-kernel-transp} holds true and take $w \in \D'\left(\TT^{N+1}\right)$ such that $\transp P w = 0$, so
\[
\left(-D_t+c(t,\xi)\right)\hat w\left(t,-\xi\right) = 0,\quad \forall \xi \in \Z^N.
\]
If $c_0(\xi) \in \C\setminus \Z$ then by~\eqref{eq:solution-ode-outside-Z}, $\hat w\left(t,-\xi\right)=0$. If $c_0(\xi) \in \Z$ then, by~\eqref{eq:solution-ode-inside-Z}, there is $C_\xi \in \C$ such that
\[
\hat w(t,\xi) = C_\xi e^{iC(t,\xi)}.
\]
Therefore
\begin{align*}
	\left\langle w, f\right\rangle &= \left(2\pi\right)^N \sum_{\xi\in\Z^N} \left\langle \hat w(t,-\xi), \hat f(t,\xi)\right\rangle\\
	&= \left(2\pi\right)^N \sum_{\xi\in\mathcal Z}C_\xi \left\langle e^{iC(t,\xi)},\hat f(t,\xi)\right\rangle\\
	&= 0,
\end{align*}
from where we conclude that $f \in \left(\ker \transp P\right)^\bot = \overline{\ran P}$.
\end{proof}

In order to deal with the quotients in~\eqref{eq:solution-ode-outside-Z} we need the next result whose proof follows the same lines of Lemma 3.1 of~\cite{bdg17}.
\begin{Lem}
	\label{Lem:DC-equivalences}
	The following conditions are equivalent:
	\begin{enumerate}
	\item There are $C>0$ and $M\in\Z_+$ such that~\eqref{eq:DC} holds true.
	\item There are $C'>0$ and $M'\in\Z_+$ such that
	\begin{align}
		\label{eq:DC-exp1}
	\left|1-e^{-2\pi i c_0(\xi)}\right| \ge C'\left(1+|\xi|\right)^{-M'},\quad \forall \xi\in\Z^N\setminus\mathcal Z.
	\end{align}
	\item There are $C''>0$ and $M''\in\Z_+$ such that
	\begin{align}
		\label{eq:DC-exp2}
	\left|e^{2\pi i c_0(\xi)}-1\right| \ge C''\left(1+|\xi|\right)^{-M''},\quad \forall \xi\in\Z^N\setminus\mathcal Z.
	\end{align}
	\end{enumerate}
\end{Lem}
\begin{proof}
	Suppose that~\eqref{eq:DC} does not hold true. Then there exist sequences $\tau_k\in\Z$ and $\xi_k \in \Z^N\setminus\mathcal Z$ such that
	\[
	\left|\tau_k + c_0\left(\xi_k\right)\right| < \left(1+\left|\xi_k\right|\right)^{-k}
	\]
	and $\xi_k\ne \xi_{k'}$ if $k\ne k'$. Indeed, we first take $C_1=\frac{1}{2}$ and $M_1=1$ to obtain $\tau_1 \in \Z$ and $\xi_1 \in \Z^N\setminus\mathcal Z$ such that
	\[
	\left|\tau_1 + c_0\left(\xi_1\right)\right| < C_1\left(1+\left|\xi_1\right|\right)^{-1} < \left(1+\left|\xi_1\right|\right)^{-1}.
	\]
	Now we take $M_2=2$ and
	\[
	C_2 = \min\left\{\frac{1}{2}, \left|\tau_1 + c_0\left(\xi_1\right)\right| \left(1+\left|\xi_1\right|\right)^{2}\right\},
	\]
	which is positive since $\xi_1 \notin \mathcal Z$. Then there exist $\tau_2\in\Z$ and $\xi_2\in\Z^N\setminus\mathcal Z$ such that
	\[
	\left|\tau_2 + c_0\left(\xi_2\right)\right| < C_2\left(1+\left|\xi_2\right|\right)^{-2} < \left(1+\left|\xi_2\right|\right)^{-2}.
	\]
	Note that $\xi_1\ne \xi_2$: otherwise,
	\begin{align*}
	\left|\tau_1-\tau_2\right| &\le \left|\tau_1+c_0\left(\xi_1\right)\right| + \left|\tau_2+c_0\left(\xi_2\right)\right| < \frac{1}{2}+\frac{1}{2} = 1,
	\end{align*}
	so $\tau_1=\tau_2$ and we obtain that
	\begin{align*}
		C_2 &\le \left|\tau_1 + c_0\left(\xi_1\right)\right| \left(1+\left|\xi_1\right|\right)^{2} = \left|\tau_2 + c_0\left(\xi_2\right)\right| \left(1+\left|\xi_2\right|\right)^{2}<C_2,
	\end{align*}
	a contradiction. Proceeding in this way we obtain the desired sequence.	
	
	Now since $\frac{1-e^{-iz}}{z}$ has a continuous extension in a neighborhood of $z=0$, there exists $c>0$ such that
	\[
	\left|1-e^{-iz}\right| \le c\left|z\right|,\quad \mbox{for $z \in \C, |z|<c^{-1}$}.
	\]
	Thus, for some $k_0\in\N$,
	\[
	\left|1-e^{-2\pi i c_0\left(\xi_k\right)}\right| = \left|1-e^{-2\pi i \left(\tau_k+c_0\left(\xi_k\right)\right)}\right| \le 2\pi c\left(1+\left|\xi_k\right|\right)^{-k}, \forall k \ge k_0
	\]
	and then there are no $C'>0$ and $M'\in\Z_+$ such that
	\[
	\left|1-e^{-2\pi i c_0\left(\xi_k\right)}\right| \ge C'\left(1+\left|\xi_k\right|\right)^{-M}, \quad \forall k \in \N.
	\]
	So~\eqref{eq:DC-exp1} implies~\eqref{eq:DC}. Conversely, if~\eqref{eq:DC-exp1} does not hold true then arguing as we have done above we conclude that there exists a sequence $\xi_k\in\Z^N\setminus\mathcal Z$ such that
	\[
	\left|1-e^{-2\pi i c_0\left(\xi_k\right)}\right| < \left(1+\left|\xi_k\right|\right)^{-k}
	\]
	and $\xi_k\ne \xi_{k'}$ if $k\ne k'$. In particular,
	\[
	 \lim _k e^{2\pi b_0\left(\xi_k\right)} = \lim_k \left| e^{-2\pi i c_0\left(\xi_k\right)} \right| = 1,
	\]
	so $b_0\left(\xi_k\right) \longrightarrow 0$. Since
	\begin{align*}
		\left|1-e^{-2\pi i c_0\left(\xi_k\right)}\right|^2 &= \left|1-e^{-2\pi i a_0\left(\xi_k\right) +2\pi b_0\left(\xi_k\right)}\right|^2\\
		&= \left(1-e^{2\pi b_0\left(\xi_k\right)}\cos\left(2\pi a_0\left(\xi_k\right)\right)\right)^2 + \left(e^{2\pi b_0\left(\xi_k\right)}\sin\left(2\pi a_0\left(\xi_k\right)\right)\right)^2\\
		&= 1 - 2e^{2\pi b_0\left(\xi_k\right)}\cos\left(2\pi a_0\left(\xi_k\right)\right) + e^{4\pi b_0\left(\xi_k\right)},
	\end{align*}
	it follows that $\cos\left(2\pi a_0\left(\xi_k\right)\right)\longrightarrow 1$. For each $k\in\N$ let $\tau_k \in \Z$ be such that
	\[
	\mathrm{dist}\left(a_0\left(\xi_k\right),\Z\right) = \left|a_0\left(\xi_k\right)+\tau_k\right|.
	\]
	We claim that
	\[
	\lim_k \tau_k + a_0\left(\xi_k\right) = 0.
	\]
	Indeed, if not there would be $\epsilon>0$ such that (after possibly considering a subsequence of $a_0\left(\xi_k\right)$)
	\[
	\epsilon \le \left|\tau_k + a_0\left(\xi_k\right)\right| \le \frac{1}{2}, \quad \forall k \in \N.
	\]
	Thus
	\[
	\epsilon 2\pi \le \left|2\pi\tau_k + 2\pi a_0\left(\xi_k\right)\right| \le \pi,
	\]
	which is not possible since $\lim_k \cos\left(2\pi\tau_k + 2\pi a_0\left(\xi_k\right)\right)=1$.
	
	Now we use that $\frac{z}{1-e^{-iz}}$ admits a continuous extension near the origin, so we can suppose that
	\[
	\left|z\right|\le c\left|1 - e^{-iz}\right|, \quad \forall |z| \le c^{-1}.
	\]
	Therefore, for $k$ large enough,
	\begin{align*}
		2\pi\left|\tau_k + c_0\left(\xi_k\right)\right| &\le c\left|1-e^{-2\pi i\left(\tau_k+c_0\left(\xi_k\right)\right)}\right|\\
		&= c\left|1-e^{-2\pi i c_0\left(\xi_k\right)}\right|\\
		&< c \left(1+\left|\xi_k\right|\right)^{-k}
	\end{align*}
	and this proves that~\eqref{eq:DC} does not hold true.
	
	The equivalence between~\eqref{eq:DC} and~\eqref{eq:DC-exp2} can be done in a similar way or one can prove directly that~\eqref{eq:DC-exp1} and~\eqref{eq:DC-exp2} are equivalent (see Remark~{\ref{Rem:CD}}).
\end{proof}
	\begin{Rem}
		\label{Rem:CD}
		Suppose that there are a number $c>0$ and a sequence $\left\{\xi_k\right\}_{k\in\N}\subset \Z^N\setminus\mathcal Z$ such that
		\begin{align}
			\label{Rem:CD-eq1}
		\left|1-e^{-2\pi i c_0\left(\xi_k\right)}\right| \le c\left(1+\left|\xi_k\right|\right)^{-k},\quad \forall k \in \N
		\end{align}
		and $\xi_k\ne \xi_{k'}$ for $k\ne k'$. Increasing $c$ and passing to a subsequence if necessary we also may suppose that
		\begin{align}
			\label{Rem:CD-eq2}
		\left|e^{2\pi i c_0\left(\xi_k\right)}-1\right| \le c\left(1+\left|\xi_k\right|\right)^{-k},\quad \forall k \in \N.
		\end{align}
		Indeed, it suffices to note that 
		\begin{align*}
			\left|e^{2\pi i c_0\left(\xi_k\right)}-1\right| & = \left|e^{2\pi i c_0\left(\xi_k\right)}\right| \left|1-e^{-2\pi i c_0\left(\xi_k\right)}\right| = \left|e^{-2\pi b_0\left(\xi_k\right)}\right| \left|1-e^{-2\pi i c_0\left(\xi_k\right)}\right| \\
			&\le e^{-2\pi b_0\left(\xi_k\right)}\left(1+\left|\xi_k\right|\right)^{-k}
		\end{align*}
		and by using that $b_0\left(\xi_k\right) \longrightarrow 0$ and by passing to a subsequence if necessary we can assume that $e^{-2\pi b_0\left(\xi_k\right)}\le 2$ for every $k\in\N$. In particular, when~\eqref{eq:DC} does not hold true then there exist $c>0$ and a sequence $\left\{\xi_k\right\}_{k\in\N}\subset\Z\setminus\mathcal Z$ such that $\xi_k\ne \xi_{k'}$ if $k\ne k'$ and~\eqref{Rem:CD-eq1} and~\eqref{Rem:CD-eq2} hold true.
	\end{Rem}

\section{Proof of Theorem~\ref{Thm:main}: the necessity of the conditions}
\label{Sec:nec}

In order to prove the necessity of conditions in Theorem~\ref{Thm:main} we recall that given a distribution $u\in\D'\left(\TT^{N+1}\right)$, $u$ is a smooth function if and only if given $\alpha\in\Z_+$ and $k\in\Z_+$ there exists $C_{\alpha,k}>0$ such that
\[
\left|D_t^\alpha\hat u\left(t,\xi\right)\right| \le C_{\alpha,k}\left(1+\left|\xi\right|\right)^{-k},\quad \forall t \in S^1, \xi \in\Z^N,
\]
where $\hat u\left(t,\xi\right)$ denotes the partial Fourier transform of $u$ with respect to the $x$-variable. We first prove the necessity of condition~\ref{Thm:main-item1} for the global solvability of $P$ (the ideas of this proof were inspired by~\cite{agkm19}):
\begin{Prop}
	\label{Prop:CD-necessaria}
	If $P$ is globally solvable then~\eqref{eq:DC} is valid.
\end{Prop}
\begin{proof}
	Suppose that~\eqref{eq:DC} does not hold true. By Remark~\ref{Rem:CD}, there exist $c>0$ and a sequence $\left\{\xi_k\right\}_{k\in\N}\subset \Z^N\setminus\mathcal Z$ such that $\xi_k\ne\xi_{k'}$ for $k\ne k'$ and
	\begin{align}\label{eq:CD-nec-eq1}
	\left|1-e^{-2\pi i c_0\left(\xi_k\right)}\right|\le c\left(1+\left|\xi_k\right|\right)^{-k}\quad \mbox{and}\quad \left|e^{2\pi i c_0\left(\xi_k\right)}-1\right|\le c\left(1+\left|\xi_k\right|\right)^{-k}, \quad \forall k \in\N.
	\end{align}
	Denote by $t_k \in [0,2\pi]$ the point where 
	\[
	t\mapsto \int_0^{t} b\left(s,\xi_k\right)\dd s
	\]
	assumes its maximum in $[0,2\pi]$. By passing to a subsequence if necessary we may assume that $t_k \longrightarrow t_0$ for some $t_0 \in [0,2\pi]$. Take $I=[\alpha,\beta]\subset [0,2\pi]$ such that $\alpha<\beta$ and $t_0 \notin I$ and by passing to a subsequence if necessary again we may assume that one of the following conditions occurs:
	\[
	t_k <\alpha, \quad \forall k \in \N
	\]
	or
	\[
	\beta<t_k,\quad \forall k \in \N.
	\]
	We fix $\varphi\in\cinfty_c((\alpha,\beta))$ such that $0\le \varphi\le 1$ and $0<\int_\alpha^{\beta}\varphi(t)\dd t$. We extend $\varphi(t)$ and $\varphi(t) C\left(t,\xi_k\right)$ as $2\pi$-periodic functions, where $C\left(t,\xi_k\right)=\int_{t_{k}}^t c\left(s,\xi_k\right)\dd s$, and we set
	\begin{align}
		\label{Prop:CD-necessaria-eq-f}
	f(t,x) = \sum_{k=1}^\infty \left(1-e^{-2\pi i c_0\left(\xi_k\right)}\right) e^{-iC\left(t,\xi_k\right)}\varphi(t)e^{ix\xi_k}
	\end{align}
	when $\beta<t_k$ for every $k\in\N$ and
	\begin{align}\label{Prop:CD-necessaria-eq-f'}
	f(t,x) = \sum_{k=1}^\infty \left(e^{2\pi i c_0\left(\xi_k\right)}-1\right) e^{-iC\left(t,\xi_k\right)}\varphi(t)e^{ix\xi_k}
	\end{align}
	when $t_k<\alpha$ for every $k\in\N$. Since $\int_{t_k}^{t} b\left(\tau,\xi_k\right)\dd\tau \le 0$ for every $t\in[0,2\pi]$ by our choice of $t_k$, it follows that
	\[
	\left|e^{-iC\left(t,\xi_k\right)}\right| = e^{\int_{t_{k}}^t b(\tau,\xi_k)}\dd\tau \le 1, \quad \forall t \in [0,2\pi],
	\]
	so, from~\eqref{eq:CD-nec-eq1}, we obtain that $f \in \cinfty\left(\TT^{N+1}\right)$. By definition, $\hat f\left(t,\xi\right)=0$ for every $\xi\in\mathcal Z$, so $f \in\overline{\ran P}$ thanks to Proposition~\ref{Prop:ker-transpose}. Since we are assuming $P$ globally solvable, there is $u\in\cinfty\left(\TT^{N+1}\right)$ such that $Pu=f$. Taking the partial Fourier transform in the $x$-variable we obtain
	\[
	\left(D_t + c\left(t,\xi_k\right)\right)\hat u\left(t,\xi_k\right) = \hat f\left(t,\xi_k\right), \quad \forall k\in\N.
	\]
	If we are in the case where $\beta<t_k$ for every $k\in\N$ then, for $s \in \left[0,2\pi\right]$, $\varphi\left(t_k-s\right)\ne 0$ implies that $t_k-s \in \left[\alpha,\beta\right]$. It follows from~\eqref{eq:solution-ode-outside-Z} that
	\begin{align*}
		\left|\hat u\left(t_k,\xi_k\right)\right| &= \left|\frac{i}{1-e^{-2\pi i c_0\left(\xi_k\right)}}\int_0^{2\pi} e^{i\left(C\left(t_k-s,\xi_k\right)-C\left(t_k,\xi_k\right)\right)} \hat f\left(t_k-s,\xi_k\right)\dd s\right|\\
		&= \left|\int_0^{2\pi} e^{iC\left(t_k-s,\xi_k\right) - i C\left(t_k,\xi_k\right)}  e^{-iC\left(t_k-s,\xi_k\right)}\varphi(t_k-s)\dd s\right|\\
		&= \int_\alpha^{\beta} \varphi(s)\dd s.
	\end{align*}
	Now if 	$t_k<\alpha$ for every $k\in\N$ then, for $s\in [0,2\pi]$, $\varphi\left(t_k+s\right) \ne 0$ implies $t_k+s \in \left[\alpha,\beta\right]$ and, by~\eqref{eq:solution-ode-outside-Z},
	\begin{align*}
		\left|\hat u\left(t_k,\xi_k\right)\right| &= \left|\frac{i}{e^{2\pi i c_0\left(\xi_k\right)}-1}\int_0^{2\pi} e^{i\left(C\left(t_k+s,\xi_k\right)-C\left(t_k,\xi_k\right)\right)} f\left(t_k+s\right)\dd s\right|\\
		&= \left|\int_0^{2\pi} e^{i\left(C\left(t_k+s,\xi_k\right)-C\left(t_k,\xi_k\right)\right)}  e^{-iC\left(t_k+s,\xi_k\right)}\varphi(t_k+s)\dd s\right|\\
		&= \int_\alpha^{\beta} \varphi(s)\dd s.
	\end{align*}
	Therefore $\int_\alpha^{\beta}\varphi(s)\dd s = \left|\hat u\left(t_k,\xi_k\right)\right|$ for every $k\in\N$, which contradicts the smoothness of $u$.
\end{proof}
In the next two results we prove the necessity of conditions $\left(\alpha\right)_D^\pm$ and $(\beta)_D$.
\begin{Prop}
	\label{Prop:alpha-nec}
	If $P$ is globally solvable then there exists $D>0$ such that given $\xi \in \Z^N\setminus \mathcal Z$, $|\xi|\ge D$, $P$ satisfies $(\alpha)^+_D$ or $(\alpha)^-_D$ at $\xi$.
\end{Prop}
\begin{proof}
Suppose that $P$ is globally solvable but given $k\in\N$, there exists $\xi_k \in \Z^N\setminus\mathcal Z$, with $\left|\xi_k\right|\ge k$, such that $P$ satisfies neither $(\alpha)_{k}^+$ nor $(\alpha)_{k}^-$ at $\xi_k$. So there are $t_{k}^+,t_{k}^- \in [0,2\pi]$ and open intervals $I_{k}^+\subset\left[t_{k}^+,t_{k}^++2\pi\right]$ and $I_{k}^-\subset\left[t_{k}^--2\pi,t_{k}^-\right]$ with $\left|\xi_k\right|^{-m}\le \left|I_{k}^\pm\right|$ and
\begin{align}
	\label{eq:negaca-psi-comp}
	B\left(t_k^{\pm},\xi_k\right) - B\left(s,\xi_k\right) > \log\left|\xi_k\right|^k, \quad \forall s \in I_k^{\pm}.
\end{align}
We may assume that $\xi_k\ne \xi_{k'}$ if $k\ne k'$ and we choose $t_k=t_k^{+}$ and $I_k = I_k^{+}$ if $0\le b_0\left(\xi_k\right)$ and $t_k=t_k^{-}$ and $I_k = I_k^{-}$ if $b_0\left(\xi_k\right)<0$. Take $J_k\subset I_k$ another open interval such that $\frac{\left|\xi_k\right|^{-m}}{2} \le \left|J_k\right|$ and $\varphi_k \in \mathscr C_c^\infty\left(I_k\right)$ such that $0\le \varphi_k\le 1$ and $\varphi_k=1$ in $J_k$. By extending $\varphi_k(t)$ and $e^{iC\left(t,\xi_k\right)}\varphi_k\left(t\right)$ as $2\pi$-periodic functions we define
\begin{align}\label{eq:Prop:alpha-nec-f}
f(t,x) \doteq \sum_{k=1}^\infty e^{i\left(C\left(t_k,\xi_k\right) -C\left(t,\xi_k\right)\right)}\varphi_k(t)e^{ix\xi_k}
\end{align}
which is smooth by~\eqref{eq:negaca-psi-comp} and belongs to $\overline{\ran P}$ thanks to Proposition~\ref{Prop:ker-transpose} since $\hat f\left(t,\xi\right)=0$ for every $\xi\in\mathcal Z$. Since we are assuming the global solvability of $P$, there exists $u\in \cinfty\left(\TT^{N+1}\right)$ such that $Pu=f$. Applying the partial Fourier transform and using~\eqref{eq:solution-ode-outside-Z} we obtain that
\begin{align*}
	\hat u\left(t,\xi_k\right) &= \frac{i}{1-e^{-2\pi i c_0\left(\xi_k\right)}} \int_0^{2\pi}e^{i\left(C\left(t-s,\xi_k\right)-C\left(t,\xi_k\right)\right)} \hat f\left(t-s,\xi_k\right) \dd s\\
	&= \frac{i}{e^{2\pi i c_0\left(\xi_k\right)}-1} \int_0^{2\pi}e^{i\left(C\left(t+s,\xi_k\right)-C\left(t,\xi_k\right)\right)} \hat f\left(t+s,\xi_k\right) \dd s.
\end{align*}
If $0\le b_0\left(\xi_k\right)$ then 
\[
\left|e^{2\pi i c_0\left(\xi_k\right)}-1\right| \le e^{-2\pi b_0\left(\xi_k\right)} + 1 \le 2.
\]
For $s\in\left[0,2\pi\right]$,  $\varphi_k\left(t_k + s\right)\ne 0$ implies that $t_k+s \in I_k=I_k^+$, so
\begin{align*}
	\left|\hat u\left(t_k,\xi_k\right)\right| &= \left|\frac{i}{e^{2\pi i c_0\left(\xi_k\right)}-1} \int_0^{2\pi}e^{i\left(C\left(t_k+s,\xi_k\right)-C\left(t_k,\xi_k\right)\right)} \hat f\left(t_k+s,\xi_k\right) \dd s\right|\\
	&= \frac{1}{\left|e^{2\pi i c_0\left(\xi_k\right)}-1\right|} \left|\int_0^{2\pi}\varphi_k(t_k+s) \dd s\right|\\
	&\ge \frac{1}{2}\left|J_k\right|\\
	&\ge \frac{\left|\xi_k\right|^{-m}}{4}.
\end{align*} 

If $b_0\left(\xi_k\right)<0$ then 
\[
\left|1-e^{-2\pi i c_0\left(\xi_k\right)}\right| \le 1 + e^{2\pi b_0\left(\xi_k\right)} \le 2
\]
 and for $s\in\left[0,2\pi\right]$,  $\varphi_k\left(t_k - s\right)\ne 0$ implies that $t_k-s \in I_k=I_k^-$, so
\begin{align*}
	\left|\hat u\left(t_k,\xi_k\right)\right| &= \frac{1}{\left|1-e^{-2\pi i c_0\left(\xi_k\right)}\right|} \left|\int_0^{2\pi}\varphi_k(t_k-s) \dd s\right|\\
	&\ge \frac{\left|\xi_k\right|^{-m}}{4}.
\end{align*}
This contradicts the smoothness of $u$ and the proof is complete.
\end{proof}

Now we proceed with the necessity of the $\left(\beta\right)_D$ condition:
\begin{Prop}
	\label{Prop:b-neceta}
	If $P$ is globally solvable then there exists $D>0$ such that given $\xi \in \mathcal Z$, $|\xi|\ge D$, $P$ satisfies $(\beta)_D$ at $\xi$.
\end{Prop}
\begin{proof}
	Suppose that $P$ is globally solvable but given $k\in\N$, there is $\xi_k\in\mathcal Z, \left|\xi_k\right|\ge k$, such that $P$ does not satisfy $(\beta)_k$ at $\xi_k$: there exist $t_k \in S^1\setminus\left\{t_{\xi_k}\right\}$ and open intervals $I_k^+\subset\left[t_k,t_{\xi_k}\right], I_k^-\subset\left[t_{\xi_k},t_k\right]$ such that $\left|I_k^{\pm}\right|\ge \left|\xi_k\right|^{-m}$ and 
	\begin{align}
	\label{eq:beta-nec-1}
	B\left(t_k,\xi_k\right) - B\left(s,\xi_k\right) > \log\left|\xi_k\right|^k,\quad \forall s \in I_k^{\pm}.
\end{align}
	We may assume that $\xi_k\ne \xi_{k'}$ for $k\ne k'$. Now we take open intervals $J_k^{\pm}\subset I_k^{\pm}$ such that $\left|J_k^{\pm}\right|\ge \frac{\left|\xi_k\right|^{-m}}{2}$ and $\varphi_k^{\pm} \in \cinfty_c\left(I_k^{\pm}\right)$ such that $0\le \varphi_k^{\pm}\le 1$, $\varphi_k^{\pm}=1$ in $J_k^\pm$ and $\int_{S^1}\varphi_k^+ = \int_{S^1}\varphi_k^-$. We set $\varphi_k = \varphi_k^+ - \varphi_k^-$ and
	\[
	f\left(t,x\right) = \sum_{k=1}^\infty e^{i\left(C\left(t_k,\xi_k\right) - C\left(t,\xi_k\right)\right)}\varphi_k(t)e^{it\xi_k},
	\]
	which is smooth thanks to~\eqref{eq:beta-nec-1}. Since $\xi_k\in\mathcal Z$ for every $k\in\N$ and
	\[
	\int_{S^1} e^{iC\left(t,\xi_k\right)} \hat f\left(t,\xi_k\right) \dd t = e^{iC\left(t_k,\xi_k\right)}\int_{S^1} \varphi_k\left(t\right) \dd t = 0,
	\]
	Proposition~\ref{Prop:ker-transpose} ensures that $f \in \overline{\ran P}$. By hypothesis, there exists $u\in \cinfty\left(\TT^{N+1}\right)$ such that $Pu=f$. We apply the partial Fourier transform and choose $t_0=t_{\xi_k}$ in~\eqref{eq:solution-ode-inside-Z} to obtain the following expression:
	\[
	\hat u\left(t,\xi_{k}\right) = C_k e^{-iC\left(t,\xi_k\right)}+i\int_{\left[t_{\xi_k},t\right]} e^{i\left(C\left(s,\xi_k\right)-C\left(t,\xi_k\right)\right)} \hat f\left(s,\xi_k\right)\dd s, \quad \forall t \in S^1, k \in \N,
	\]
	for some $C_k \in \C$. Taking $t=t_{\xi_k}$,
	\[
	C_k = \hat u\left(t_{\xi_k},\xi_k\right) e^{iC\left(t_{\xi_k},\xi_k\right)},
	\]
	so
	\begin{align*}
		\left|\int_{\left[t_{\xi_k},t_k\right]} e^{i\left(C\left(s,\xi_k\right)-C\left(t_k,\xi_k\right)\right)} \hat f\left(s,\xi_k\right)\dd s\right| &=	\left|\hat u\left(t_k,\xi_{k}\right) - \hat u\left(t_{\xi_k},\xi_k\right) e^{i\left(C\left(t_{\xi_k},\xi_k\right)-C\left(t_k,\xi_k\right)\right)}\right|\\
		&\le \left|\hat u\left(t_k,\xi_k\right)\right| + \left|\hat u\left(t_{\xi_k},\xi_k\right)\right| e^{-B\left(t_{\xi_k},\xi_k\right) + B\left(t_k,\xi_k\right)}\\
		&\le  \left|\hat u\left(t_k,\xi_k\right)\right| + \left|\hat u\left(t_{\xi_k},\xi_k\right)\right|, \quad \forall k \in \N.
	\end{align*}
	The smoothness of $u$ ensures that given $n \in \Z_+$ there exists $C_n>0$ such that
	\[
	\left|\int_{\left[t_{\xi_k},t_k\right]} e^{i\left(C\left(s,\xi_k\right)-C\left(t_k,\xi_k\right)\right)} \hat f\left(s,\xi_k\right)\dd s\right| \le C_n\left(1+\left|\xi_k\right|\right)^{-n},\quad \forall k \in \N.
	\]
	But note that
	\begin{align*}
		\left|\int_{\left[t_{\xi_k},t_k\right]} e^{i\left(C\left(s,\xi_k\right)-C\left(t_k,\xi_k\right)\right)} \hat f\left(s,\xi_k\right)\dd s\right| &= \left|\int_{\left[t_{\xi_k},t_k\right]} \varphi_k\left(s\right)\dd s\right|\\
		&= \left|\int_{\left[t_{\xi_k},t_k\right]} \varphi_k^-\left(s\right)\dd s\right|\\
		&\ge \frac{\left|\xi_k\right|^{-m}}{2},
	\end{align*}
	which gives a contradiction.
\end{proof}

\section{Proof of Theorem~\ref{Thm:main}: the sufficiency of the conditions}
\label{Sec:suf}

Now we assume the conditions~\ref{Thm:main-item1},~\ref{Thm:main-item2} and~\ref{Thm:main-item3} of Theorem~\ref{Thm:main} and we prove that $P$ is globally solvable. We start with the following
\begin{Prop}
	\label{Prop:decaimento-sem-derivadas}
	Suppose that $f \in \cinfty\left(\TT^{N+1}\right)$ and $\left\{u_\xi\right\}_{\xi\in\Z^{N}}\subset \cinfty\left(S^1\right)$ is a family of functions such that
	\[
	\left(D_t + c(t,\xi)\right) u_\xi(t) = \hat f(t,\xi), \quad \forall t \in S^1,\xi\in\Z^{N}.
	\]
	Then the sum
	\[
	\sum_{\xi\in\Z^N} u_\xi(t)\otimes e^{ix\xi}
	\]
	belongs to $\cinfty\left(\TT^{N+1}\right)$ if and only if given $k\in\Z_+$, there exists $C_k>0$ such that
	\[
	\left|u_\xi(t)\right|\le C_k\left(1+|\xi|\right)^{-k},\quad \forall t \in S^1, \xi \in \Z^N.
	\]
\end{Prop}
\begin{proof}
	We need to prove that given $\alpha \in \Z_+$ and $k\in\Z_+$ there exists $C_{\alpha,k}>0$ such that
	\[
	\left|D_t^\alpha u_\xi\left(t\right)\right| \le C_{\alpha,k}\left(1+\left|\xi\right|\right)^{-k},\quad \forall t \in S^1, \xi \in \Z^N. 
	\]
	This can be done by induction in $\alpha$. For $\alpha=0$ we consider $C_{\alpha,0}=C_k$. Suppose that the result holds true for every $0\le \beta \le \alpha$ and let us prove that it is also valid for $\alpha+1$. Since 
	\[
	D_t u_\xi\left(t\right) = - c\left(t,\xi\right)u_\xi\left(t\right)+\hat f\left(t,\xi\right),
	\]
	we have that
	\[
	\left|D_t^{\alpha+1}u_\xi\left(t\right)\right| \le \sum_{\beta\le\alpha}\binom{\alpha}{\beta} \left|D_t^{\alpha-\beta}c\left(t,\xi\right)\right|\left|D_t^\beta u_\xi\left(t\right)\right| + \left|D_t^\alpha \hat f\left(t,\xi\right)\right|, \quad \forall t \in S^1, \xi \in \Z^N. 
	\]
	Given $k\in\Z_+$, there exists $D_{\alpha,k}>0$ such that 
	\[
	\left|D_t^\alpha \hat f\left(t,\xi\right)\right| \le D_{\alpha,k}\left(1+\left|\xi\right|\right)^{-k},\quad \forall t \in S^1, \xi \in \Z^N
	\]
	and there is $D>0$ such that
	\[
	\left|D_t^\gamma c\left(t,\xi\right)\right| \le D\left(1+\left|\xi\right|\right)^m,\quad \forall t \in S^1, \xi \in \Z^N, 0\le \gamma\le \alpha.
	\]
	Therefore, by the inductive hypothesis,
	\begin{align*}
		\left|D_t^{\alpha+1}u_\xi\left(t\right)\right| &\le \sum_{\beta\le\alpha}\binom{\alpha}{\beta} D\left(1+\left|\xi\right|\right)^{m} C_{\beta,k+m}\left(1+\left|\xi\right|\right)^{-k-m} + D_{\alpha,k}\left(1+\left|\xi\right|\right)^{-k}\\
		&\le D'\left(1+\left|\xi\right|\right)^{-k},\quad \forall t \in S^1, \xi\in\Z^N,
	\end{align*}
	where $D'$ depends only on $\alpha$ and $k$.
\end{proof}

Given $X\subset \Z^N$, we denote by $\cinfty_X\left(\TT^{N+1}\right)$ the subspace of $\cinfty\left(\TT^{N+1}\right)$ formed by the functions $f$ such that $\hat f(t,\xi) = 0$ if $\xi\in\Z^N\setminus X$. We can decompose
\[
\cinfty\left(\TT^{N+1}\right) = \cinfty_\mathcal Z\left(\TT^{N+1}\right) \oplus \cinfty_{\Z^{N}\setminus \mathcal Z}\left(\TT^{N+1}\right)
\]
and, in order to prove that $P$ is globally solvable, it suffices to prove that for each $f\in\overline{\ran P}\cap \cinfty_X\left(\TT^{N+1}\right)$, where $X=\mathcal Z$ or $X=\Z^N\setminus \mathcal Z$, there exists $u\in\cinfty_X\left(\TT^{N}\right)$ such that $Pu=f$. In this case we say that $P$ is globally solvable in $\cinfty_X\left(\TT^{N+1}\right)$ -- also note that each $\cinfty_X\left(\TT^{N+1}\right)$ is invariant by $P$.

What we will actually prove is that the conditions~\ref{Thm:main-item1} and \ref{Thm:main-item2} are equivalent to the global solvability of $P$ in $\cinfty_{\Z^N\setminus\mathcal Z}\left(\TT^N\right)$ and the condition~\ref{Thm:main-item3} is equivalent to the global solvability of $P$ in $\cinfty_{\mathcal Z}\left(\TT^{N+1}\right)$. The reader should note that the functions~\eqref{Prop:CD-necessaria-eq-f},~\eqref{Prop:CD-necessaria-eq-f'} and~\eqref{eq:Prop:alpha-nec-f} belong to $\cinfty_{\Z^N\setminus\mathcal Z}\left(\TT^N\right)$ and in the proof of Proposition~\ref{Prop:b-neceta} the function $f$ belongs to $\cinfty_{\mathcal Z}\left(\TT^{N+1}\right)$.

{\bf Solvability in $\cinfty_{\Z^N\setminus\mathcal Z}\left(\TT^{N+1}\right)$}: suppose that~\eqref{eq:DC} holds true and that there exists $D>0$ such that for each $\xi\in\Z^N\setminus\mathcal Z$, with $\left|\xi\right|\ge D$, $P$ satisfies $(\alpha)_D^+$ or $(\alpha)_D^-$ at $\xi$. Let $f\in \cinfty_{\Z^N\setminus\mathcal Z}\left(\TT^{N+1}\right)\subset \overline{\ran P}$ and we will prove that there is $u\in\cinfty_{\Z^N\setminus\mathcal Z}\left(\TT^{N+1}\right)$ such that $Pu=f$. For each $\xi\in\Z^N\setminus\mathcal Z$, since $c_0\left(\xi\right) \notin \Z$, we apply the partial Fourier transform in the $x$-variable in the equation $Pu=f$ and~\eqref{eq:solution-ode-outside-Z} gives that $\hat u(t,\xi)$ must satisfy
\begin{align*}
	\hat u\left(t,\xi\right) &= \frac{i}{1-e^{-2\pi i c_0\left(\xi\right)}} \int_0^{2\pi}e^{i\left(C\left(t-s,\xi\right)-C\left(t,\xi\right)\right)} \hat f\left(t-s,\xi\right) \dd s\\
	&= \frac{i}{e^{2\pi i c_0\left(\xi\right)}-1} \int_0^{2\pi}e^{i\left(C\left(t+s,\xi\right)-C\left(t,\xi\right)\right)} \hat f\left(t+s,\xi\right) \dd s.
\end{align*}
Since $\hat f(t,\xi)=0$ for $\xi\in\mathcal Z$, we also can take $\hat u(t,\xi) = 0$ and we only have to prove that
\[
u(t,x) \doteq \sum_{\xi\in\Z^N} \hat u\left(t,\xi\right) \otimes e^{ix\xi}
\]
belongs to $\cinfty\left(\TT^{N+1}\right)$. From Proposition~\ref{Prop:decaimento-sem-derivadas}, the hypothesis~\eqref{eq:DC} and Lemma~\ref{Lem:DC-equivalences}, it suffices to prove that given $k\in\Z_+$, there exists $C_k>0$ such that for each $\xi \in \Z^N\setminus \mathcal Z, |\xi|\ge D$,
\[
\left|\int_0^{2\pi}e^{i\left(C\left(t-s,\xi\right)-C\left(t,\xi\right)\right)} \hat f\left(t-s,\xi\right) \dd s\right| \le C_k\left(1+|\xi|\right)^{-k},\quad \forall t \in [0,2\pi]
\]
or
\[
\left|\int_0^{2\pi}e^{i\left(C\left(t+s,\xi\right)-C\left(t,\xi\right)\right)} \hat f\left(t+s,\xi\right) \dd s\right| \le C_k\left(1+|\xi|\right)^{-k},\quad \forall t \in [0,2\pi].
\]
We know that there is $D_k>0$ such that 
\begin{align}
	\label{eq:suf-est-f}
\left|\hat f\left(t,\xi\right)\right| \le D_k\left(1+|\xi|\right)^{-k-D},\quad \forall t \in [0,2\pi], \xi \in \Z^N\setminus\mathcal Z.
\end{align}
We claim that we can take $C_k=2\pi e^{c} D_k$, where $c>0$ satisfy
\begin{align}
	\label{eq:suf-est-c}
\left|b\left(t,\xi\right)\right|\le c\left|\xi\right|^m,\quad \forall t \in [0,2\pi], \xi\ne 0.
\end{align}
 Indeed, fix $\xi \in \Z^N\setminus\mathcal Z$, $\left|\xi\right|\ge D$ and suppose that $(\alpha)_D^+$ holds true at $\xi$. We shall analyze the term
\[
(\star)=\int_0^{2\pi} e^{i\left(C\left(t+s,\xi\right)-C\left(t,\xi\right)\right)} \hat f\left(t+s,\xi\right) \dd s
\]
for each $t \in [0,2\pi]$ fixed. Take $n\in\Z_+$ the largest integer such that $n\left|\xi\right|^{-m} < 2\pi$, define $\alpha=2\pi - n\left|\xi\right|^{-m}$ and we can write
\begin{align*}
	(\star )= \underbrace{\int_0^\alpha e^{i\left(C\left(t+s,\xi\right)-C\left(t,\xi\right)\right)} \hat f\left(t+s,\xi\right) \dd s}_{\doteq (\star)_0}  + \sum_{j=1}^{n} \underbrace{\int_{(j-1)\left|\xi\right|^{-m}+\alpha}^{j\left|\xi\right|^{-m}+\alpha} e^{i\left(C\left(t+s,\xi\right)-C\left(t,\xi\right)\right)} \hat f\left(t+s,\xi\right) \dd s}_{\doteq(\star)_j}.
\end{align*}
Since $2\pi \le (n+1)\left|\xi\right|^{-m}$ it follows that $\alpha\le \left|\xi\right|^{-m}$, so
\begin{align*}
	\left|(\star)_0\right| &= \left|\int_0^\alpha e^{i\int_t^{t+s} c\left(\tau,\xi\right)\dd \tau} \hat f\left(t+s,\xi\right) \dd s\right|\\
	&\le D_k\left(1+|\xi|\right)^{-k-D} \int_0^{\alpha} e^{-\int_t^{t+s} b\left(\tau,\xi\right)\dd \tau}\dd s\\
	&\le D_k\left(1+|\xi|\right)^{-k} \int_0^{\alpha} e^{sc\left|\xi\right|^m}\dd s\\
	&\le  D_k\left(1+|\xi|\right)^{-k} \int_0^{\alpha} e^{\alpha c\left|\xi\right|^m}\dd s\\
	&\le D_k\left(1+|\xi|\right)^{-k} \int_0^{\alpha} e^{c}\dd s\\
	&\le \alpha e^{c}D_k\left(1+|\xi|\right)^{-k}.
\end{align*}
Now for each $j\in\{1,\ldots,n\}$, since the length of the interval 
\[
I_{j,\xi} = \left(t+\left(j-1\right)\left|\xi\right|^{-m}+\alpha, t+j\left|\xi\right|^{-m}+\alpha\right) \subset \left[t,t+2\pi\right]
\]
is $\left|\xi\right|^{-m}$, it follows by hypothesis that there is $s_j \in I_{j,\xi}$ such that
\[
B\left(t,\xi\right)-B\left(s_j,\xi\right) \le D\log|\xi|.
\]
We define $s_j' = s_j-t \in \left(\left(j-1\right)\left|\xi\right|^{-m}+\alpha, j\left|\xi\right|^{-m}+\alpha\right)$ and then
\begin{align*}
	\left|(\star)_j\right| &\le \int_{(j-1)\left|\xi\right|^{-m}+\alpha}^{j\left|\xi\right|^{-m}+\alpha} e^{B\left(t,\xi\right)-B\left(s_j,\xi\right)}  e^{B\left(s_j,\xi\right)-B\left(t+s,\xi\right)} \left|\hat f\left(t+s,\xi\right)\right| \dd s\\
	&\le D_k\left(1+|\xi|\right)^{-k-D}\int_{(j-1)\left|\xi\right|^{-m}+\alpha}^{j\left|\xi\right|^{-m}+\alpha} |\xi|^{D}  e^{\int_{t+s}^{t+s_j'} b\left(\tau,\xi\right)\dd \tau}  \dd s\\
	&\le D_k\left(1+|\xi|\right)^{-k} \int_{(j-1)\left|\xi\right|^{-m}+\alpha}^{j\left|\xi\right|^{-m}+\alpha} e^{|s_j'-s|c|\xi|^m}\dd s\\
	&\le e^{c}D_k\left(1+|\xi|\right)^{-k} |\xi|^{-m}.
\end{align*}
Therefore
\begin{align*}
	\left|(\star)\right|&\le \alpha e^{c} D_k\left(1+|\xi|\right)^{-k} + \sum_{j=1}^n e^{c}D_k\left(1+|\xi|\right)^{-k} |\xi|^{-m}\\
	&= e^{c} D_k\left(1+|\xi|\right)^{-k} \left(\alpha+ n|\xi|^{-m}\right)\\
	&= 2\pi e^{c} D_k\left(1+|\xi|\right)^{-k}.
\end{align*}
When $(\alpha)_D^-$ holds true at $\xi$ we work in a similar way with the term
\[
\int_0^{2\pi} e^{i\left(C\left(t-s,\xi\right)-C\left(t,\xi\right)\right)} \hat f\left(t-s,\xi\right) \dd s.
\]

{\bf Solvability in $\cinfty_{\mathcal Z}\left(\TT^{N+1}\right)$}: suppose that there exists $D>0$ such that for every $\xi \in \mathcal Z, |\xi|\ge D$, $P$ satisfies $(\beta)_{D}$ at $\xi$ and let $f \in \overline{\ran P}\cap\cinfty_\mathcal Z\left(\TT^{N+1}\right)$. 
In this case we define $\hat u\left(t,\xi\right) = 0$ for $\xi\in\Z^N\setminus\mathcal Z$ and
\[
\hat u(t,\xi) \doteq i\int_{\left[t_{\xi},t\right]} e^{i\left(C(s,\xi)-C(t,\xi)\right)}\hat f(s,\xi)\dd s, \quad t \in S^1, \xi\in\mathcal Z.
\]
Thanks to Proposition~\ref{Prop:ker-transpose}, $\hat u\left(t,\xi\right)\in\cinfty\left(S^1\right)$ for every $\xi\in\mathcal Z$. Moreover,~\eqref{eq:solution-ode-inside-Z} ensures that $\left(D_t + c\left(t,\xi\right)\right)\hat u\left(t,\xi\right)=\hat f\left(t,\xi \right)$ for every $\xi\in\Z^N$, so in order to show that 
\[
u=\sum_{\xi\in\mathcal Z} \hat u\left(t,\xi\right) \otimes e^{ix\xi}
\]
 defines a solution for $Pu=f$ we will prove (see Proposition~\ref{Prop:decaimento-sem-derivadas}) that 
\begin{align}
	\label{eq:ineq-u-smooth}
\left|\hat u(t,\xi)\right| \le 2\pi e^{c} D_k\left(1+|\xi|\right)^{-k}, \quad \forall \xi \in \mathcal Z, \left|\xi\right|\ge D, t \in S^1, k \in \Z_+,
\end{align}
where $D_k$ and $c$ satisfy~\eqref{eq:suf-est-f} and~\eqref{eq:suf-est-c} respectively. By the definition of $\hat u(t,\xi)$, the inequality~\eqref{eq:ineq-u-smooth} is trivially satisfied for $t=t_{\xi}$, so we will work only with $t \in S^1\setminus\{t_{\xi}\}$.
Fix $\xi\in\mathcal Z$, $|\xi|\ge D$, $t\in S^1\setminus\left\{t_\xi\right\}$ and $k\in\Z_+$ and recall that we are assuming that $P$ satisfies $\left(\beta\right)_D$ at $\xi$. Suppose that given an open interval $I^-\subset\left[t_{\xi},t\right]$, with $\left|I^-\right| \ge \left|\xi\right|^{-m}$, there is $s \in I^-$ such that
\[
B\left(t,\xi\right) - B\left(s,\xi\right) \le D\log\left|\xi\right|.
\]
Note that
\begin{align*}
	\left|\hat u(t,\xi)\right| &\le D_k\left(1+\left|\xi\right|\right)^{-k-D} \left|\left[t_\xi,t\right]\right| \sup_{s \in \left[t_\xi,t\right]} e^{B\left(t,\xi\right)-B\left(s,\xi\right)} \\
	&\le D_k\left(1+\left|\xi\right|\right)^{-k} \left|\left[t_\xi,t\right]\right| \sup_{s \in \left[t_\xi,t\right]} e^{\int_{\left[s,t\right]} b\left(\tau,\xi\right)\dd\tau}\\
	&\le D_k\left(1+\left|\xi\right|\right)^{-k} \left|\left[t_\xi,t\right]\right| \sup_{s \in \left[t_\xi,t\right]} e^{\left|\left[s,t\right]\right| c \left|\xi\right|^{m}}.
\end{align*}
In particular, if $\left|\left[t_\xi,t\right]\right| \le \left|\xi\right|^{-m}$ then
\[
\left|\hat u(t,\xi)\right| \le 2\pi D_k \left(1+\left|\xi\right|\right)^{-k}  e^{c}.
\]
If $\left|\left[t_\xi,t\right]\right| > \left|\xi\right|^{-m}$ we proceed as we have done for the solvability in $\cinfty_{\Z^N\setminus\mathcal Z}\left(\TT^{N+1}\right)$: we first decompose
\[
\left[t_\xi,t\right] = \left[t_\xi,\alpha\right] \cup I_1\cup\ldots\cup I_n,
\]
where $I_1,\ldots,I_n$ are consecutive closed intervals (in the anticlockwise orientation of $S^1$) with disjoint interiors, $\left|\left[t_\xi,\alpha\right] \right| \le \left|\xi\right|^{-m}$ and $\left|I_j\right|=|\xi|^{-m}$. By hypothesis, there is $s_j \in I_j$ such that
\[
B\left(t,\xi\right) - B\left(s_j,\xi\right) \le D\log\left|\xi\right|,\quad \forall j \in \{1,\ldots,n\}. 
\]
Hence, by the previous analysis,
\begin{align*}
	\left|\hat u(t,\xi)\right| &\le \left|\int_{\left[t_\xi,\alpha\right]}e^{i\left(C(s,\xi)-C(t,\xi)\right)}\hat f(s,\xi)\dd s\right|+ \sum_{j=1}^n \left|\int_{I_j} e^{i\left(C(s,\xi)-C(t,\xi)\right)}\hat f(s,\xi)\dd s\right|\\
	 &\le D_k\left(1+\left|\xi\right|\right)^{-k} \left|\left[t_\xi,\alpha\right]\right|e^{c}+ \sum_{j=1}^n \underbrace{\left|\int_{I_j} e^{i\left(C(s,\xi)-C(t,\xi)\right)}\hat f(s,\xi)\dd s\right|}_{\doteq (\dagger_j)}.
\end{align*}
Note that
\begin{align*}
	(\dagger_j) &\le \left|I_j\right| \sup_{s \in I_j} \left|e^{i\left(C(s,\xi)-C\left(s_j,\xi\right) + C\left(s_j,\xi\right)-C(t,\xi)\right)}\hat f(s,\xi)\right|\\
	&\le D_k\left(1+\left|\xi\right|\right)^{-k-D} \left|I_j\right|  \sup_{s\in I_j} e^{B\left(s_j,\xi\right)-B(s,\xi)} e^{B\left(t,\xi\right)-B\left(s_j,\xi\right)}\\
	&\le D_k\left(1+\left|\xi\right|\right)^{-k} \left|I_j\right|\sup_{s\in I_j} e^{\int_{[s,s_j]} b\left(\tau,\xi\right)\dd\tau} \\
	&\le D_k\left(1+\left|\xi\right|\right)^{-k} \left|I_j\right| \sup_{s \in I_j} e^{\left|\left[s,s_j\right]\right|c\left|\xi\right|^{m}} \\ 
	&= e^{c}D_k\left(1+\left|\xi\right|\right)^{-k} \left|I_j\right|.
\end{align*}
Therefore
\begin{align*}
	\left|\hat u\left(t,\xi\right)\right| &\le D_k\left(1+\left|\xi\right|\right)^{-k} \left|\left[t_\xi,\alpha\right]\right|e^{c}+ \sum_{j=1}^n e^{c}D_k\left(1+\left|\xi\right|\right)^{-k+D} \left|I_j\right|\\
	&= e^{c} D_k\left(1+\left|\xi\right|\right)^{-k}\left|\left[t_\xi,t\right]\right|\\
	&\le 2\pi e^{c} D_k \left(1+\left|\xi\right|\right)^{-k}.
\end{align*}
The case where given $I^+\subset\left[t,t_{\xi}\right]$, with $\left|I^+\right| \ge \left|\xi\right|^m$, there is $s \in I^+$ such that
\[
B\left(t,\xi\right) - B\left(s,\xi\right) \le D\log\left|\xi\right|
\]
can be done analogously since 
\[
\hat u\left(t,\xi\right) = - i \int_{\left[t,t_\xi\right]} e^{i\left(C\left(s,\xi\right)-C\left(t,\xi\right)\right)}\hat f\left(s,\xi\right)\dd s
\]
thanks to Proposition~\ref{Prop:ker-transpose}.

\section{Applications and examples}
\label{Sec:app}
In this section we present applications of Theorem~\ref{Thm:main}.

%

\subsection{Operators with negative order} 

Suppose that $b(t,D_x)$ has order $m<0$. Then $2\pi < |\xi|^{-m}$ for every $\xi\in\Z^N$ except a finite set. So $(\alpha)^{\pm}_1$ and $(\beta)_1$ are trivial in this case and $P$ is globally solvable if and only if the condition~\eqref{eq:DC} holds true. Alternatively, there are $C,C'>0$ such that
\[
B(t,\xi)-B(s,\xi) \le C(1+|\xi|)^m <  C'\log|\xi|, \quad \forall t,s \in [0,2\pi], \left|\xi\right| > 1.
\]

\subsection{Homogeneous operators}

In this section we work with positively homogenegous operators of degree $m$, that is $c\left(t,D_x\right)$ has order $m$ and
\begin{align}\label{eq:homog}
c(t,n\xi) = n^m c\left(t,\xi\right), \quad \forall t \in S^1, \xi \in \Z^N, n \in\N.
\end{align}
Note that $c\left(t,D_x\right)$ is positively homogeneous of degree $m$ if and only if its real and imaginary part are both positively homogeneous of degree $m$.

	\begin{Prop}
		\label{Prop:homogeneous-complementar-Z-necessary}
		Let $c\left(t,D_x\right)$ be a positively homogeneous pseudodifferential operator of order $m\in\N$ such that there exists $\xi_0\in\Z^N\setminus\{0\}$ satisfying
		\begin{enumerate}
			\item $c_{0}\left(\xi_0\right) \notin \Z$;
			\item $b\left(t,\xi_0\right)$ changes signal.
		\end{enumerate}
		Then $P=D_t+c\left(t,D_x\right)$ is not globally solvable.
	\end{Prop}
	\begin{proof}
		Since $c_0\left(\xi_0\right)\notin \Z$ and
		\[
		c_0\left(n\xi_0\right) = \left(2\pi\right)^{-1}\int_0^{2\pi} c\left(t,n\xi_0\right)\dd t = n^m c_0\left(\xi_0\right),\quad \forall n \in\N,
		\]
		it follows that there exist an infinite set $\mathcal N\subset \N$ such that $c_0\left(n\xi_0\right) \notin \Z$ for every $n\in\mathcal N$.
		
		Suppose that there is $D>0$ such that $P$ satisfies $(\alpha)_D^+$ or $(\alpha)_D^-$ at $n\xi_0$ for each $n\in\mathcal N$ such that $n\left|\xi_0\right|\ge D$. Let us say that $P$ satisfies $(\alpha)_D^+$ at $n\xi_0$ for an infinite number of $n\in\mathcal N$ such that $n\left|\xi_0\right|\ge D$. Changing $\mathcal N$ if necessary, we can suppose that $P$ satisfies $(\alpha)_D^+$ for every $n\in\mathcal N$. 
		
		By hypothesis, $b\left(t,\xi_0\right)$ changes signal, so there exists $t_0$ in $[0,2\pi]$ such that
				\[
				b\left(t_0,\xi_0\right)<0.
				\]
		If we consider the primitives
			\[
			B\left(t,\xi\right) = \int_{t_0}^t b\left(\tau,\xi\right)\dd \tau, \quad \xi\in\Z^N,
			\]	
		then for any open interval $I^+\subset\left[t_0,t_0+2\pi\right]$ such that $\left|I^+\right|\ge n^{-m}\left|\xi_0\right|^{-m}$, there exists $s_n \in I^+$ such that
		\begin{align*}
			-n^{m}B\left(s_n,\xi_0\right) = -B\left(s_n,n\xi_0\right)= B\left(t_0,n\xi_0\right) - B\left(s_n,n\xi_0\right) \le D \log\left|n\xi_0\right|,
		\end{align*}
		that is
			\begin{align}
			\label{eq:alpha-cond-homog-1'}
			-B\left(s_n,\xi_0\right)\le D n^{-m}\log\left|n\xi_0\right|.
		\end{align}
		Since $B(t_0,\xi_0)=0$ and $\partial_{t}B(t_0,\xi_0)<0$, there exist $0<\delta<2\pi$ such that
		\[
		B\left(t,\xi_0\right)<0, \quad \forall t \in \left(t_0,t_0+\delta\right).
		\]
		Decreasing $\delta>0$ if necessary we may assume that
		\[
		0<\theta\doteq \min_{s \in \left[t_0+\frac{\delta}{2},t_0+\delta\right]} -B\left(s,\xi_0\right).
		\]
		Discarding a finite number of terms in $\mathcal N$ we also can assume that
		\[
		n^{-m}\left|\xi_0\right|^{-m} \le \frac{\delta}{2}, \quad \forall n\in\mathcal N.
		\]
		It follows from \eqref{eq:alpha-cond-homog-1'} that for each $n\in\mathcal N$, there exists $s_n\in\left(t_0+\frac{\delta}{2},t_0+\delta\right)$ such that
		\[
		\theta \le -B\left(s_n,\xi_0\right)\le n^{-m} D \log\left|n\xi_0\right|.
		\]
		This gives a contradiction since $\mathcal N$ is infinite and $0<m$.
		
		If $(\alpha)_D^-$ is satisfied for an infinite set $\mathcal N$, we work in a similar way with a point $t_1\in[0,2\pi]$ where $0<b\left(t_1,\xi_0\right)$.
	\end{proof}
	Proposition~\ref{Prop:homogeneous-complementar-Z-necessary} can be improved in several ways. For instance, we actually have proved that when $b(\cdot,\xi_0)$ assumes a negative value then there is no $D>0$ such that $P$ satisfies $\left(\alpha\right)_D^+$ at $n\xi_0$ for infinitely many $n\in\N$ such that $n\xi_0 \notin \mathcal Z$. Analogously, if $b\left(\cdot,\xi_0\right)$ assumes a positive value then there is no $D>0$ such that $\left(\alpha\right)_D^-$ holds true at $n\xi_0$ for infinitely many $n\in\N$ such that $n\xi_0 \notin\mathcal Z$. Another improvement is a necessary condition for the global solvability of $P$ given in terms of the principal part of $b\left(t,D_x\right)$:
		\begin{Prop}
				\label{Prop:homogeneous-complementar-Z-necessary-imag}
				Let $b\left(t,D_x\right)$ be a pseudodifferential of order $m>0$ such that $b\left(t,\xi\right)=b_m\left(t,\xi\right)+b'\left(t,\xi\right)$, where $b_m$ is positively homogeneous of degree $m$ and $b'(x,\xi)$ is a symbol of order $m'<m$. Suppose that there exists $\xi_0 \in \Z^N\setminus\{0\}$ such that
				\begin{enumerate}
						\item $b_{m0}\left(\xi_0\right) = (2\pi)^{-1}\int_0^{2\pi} b_m\left(t,\xi_0\right)\dd t \ne 0$;
						\item $b_m\left(t,\xi_0\right)$ changes signal.
					\end{enumerate}
				Then $P=D_t+c\left(t,D_x\right)$ is not globally solvable.
			\end{Prop}
		Indeed, we first take $c>0$ such that
				\[
				\left|b'\left(t,n\xi_0\right)\right| \le cn^{m'}, \quad \forall t \in [0,2\pi],n \in \N.				\]
		Since $m'<m$, $b_m$ is positively homogeneous of degree $m$ and $b_{m0}\left(\xi_0\right)\ne 0$, there exists $n_0 \in \N$ such that $b_0\left(n\xi_0\right) \ne 0$ for every $n\ge n_0$. This implies that $n\xi_0 \notin\mathcal Z$ for $n\ge n_0$. By hypothesis $b_m\left(t,\xi_0\right)$ changes signal, so there exist $t_0,t_1$ in $[0,2\pi]$ such that
				\[
				b_m\left(t_0,\xi_0\right)<0<b_m\left(t_1,\xi_0\right).
				\] 
				As before first assume that there exist an infinite set $\mathcal N \subset \N$ such that, for some $D>0$, $P$ satisfies $(\alpha)_D^+$ at $n\xi_0$ and $n\xi_0 \notin \mathcal Z$ for every $n\in\mathcal N$. We consider
				\[
				B_m\left(t,n\xi_0\right) = \int_{t_0}^t b_m\left(\tau,n\xi_0\right)\dd \tau \quad\mbox{and}\quad B'\left(t,n\xi_0\right) = \int_{t_0}^t b'\left(\tau,n\xi_0\right)\dd\tau, \quad n \in \N,
				\]
				and since $B_m(t_0,\xi_0)=0$ and $\partial_{t}B_m(t_0,\xi_0)<0$, there exists $0<\delta<2\pi$ such that
				%
				%
				\[
				0<\theta\doteq \min_{s \in \left[t_0+\frac{\delta}{2},t_0+\delta\right]} -B_m\left(s,\xi_0\right).
				\]
				Thus, for $s \in \left[t_0+\frac{\delta}{2},t_0+\delta\right]$,
				\begin{align*}
						-B\left(s,n\xi_0\right) &= -B_m\left(s,n\xi_0\right) - B'\left(s,n\xi_0\right)\\
						&\ge n^m\theta - 2\pi cn^{m'}\\
						&= n^m\left(\theta-2\pi cn^{m'-m}\right).
					\end{align*}
				Changing $\mathcal N$ if necessary we may assume that
				\[
				-B\left(s,n\xi_0\right) \ge \frac{\theta}{2}n^m,\quad \forall n\in \mathcal N, s \in \Big[t_0+\frac{\delta}{2},t_0+\delta\Big]
				\]
				and also that
				\[
				n^{-m}\left|\xi_0\right|^{-m} \le \frac{\delta}{2}, \quad \forall n \in\mathcal N,
				\]
				which gives a contradiction as in the proof of Proposition~\ref{Prop:homogeneous-complementar-Z-necessary}.

\begin{Prop}
	\label{Prop:homogeneous-Z-necessary}
	Let $c\left(t,D_x\right)$ be a positively homogeneous pseudodifferential operator of order $m\in\N$ such that $c_0\left(\xi_0\right) \in \Z$ for some $\xi_0 \in\Z^N\setminus\{0\}$. If there exists $\lambda\in\R$ such that
	\[
	\Omega_\lambda \doteq\left\{t\in S^1 : B\left(t,\xi_0\right) < \lambda\right\}
	\]
	is not connected then $P=D_t+c\left(t,D_x\right)$ is not globally solvable. 
\end{Prop}
\begin{proof}
	Note that
	\begin{align*}
		c_{0}\left(n\xi_0\right) &= (2\pi)^{-1}\int_0^{2\pi} c_0\left(t,n\xi_0\right)\dd t = n^{m} c_0\left(\xi_0\right) \in \Z,\quad \forall n \in \N,
	\end{align*}
	so $n\xi_0 \in \mathcal Z$ for every $n \in \N$.
	
	By hypothesis there are two nonempty open components $I^+,I^- \subset S^1$ of $\Omega_\lambda$. Then $I^+$ and $I^-$ are open and disjoint intervals in $S^1$ and $B\left(t,\xi_0\right)<\lambda$ for every $t \in I^+\cup I^-$. Since $\Omega_\lambda \ne S^1$ there is also a point $t_0 \in S^1$ such that $\lambda\le B\left(t_0,\xi_0\right)$. Recall that $t_{\xi_0}\in S^1$ denotes a point where $B\left(\cdot,\xi_0\right)$ assumes its maximum in $S^1$, so $t_{\xi_0}\notin I^+\cup I^-$. The set $S^1\setminus \left(I^+\cup I^-\right)$ has two components and one of them contains $t_{\xi_0}$. We can assume that $t_0$ belongs to the other component and to $\partial I^+$. Thus if we fix $s^\pm \in I^\pm$, we obtain that
	\[
	\lambda = B\left(t_0,\xi_0\right) > B\left(s^\pm,\xi_0\right).
	\]
	Now we take $\lambda'\in \R$ such that
	\[
	\lambda > \lambda' > B\left(s^\pm,\xi_0\right)
	\]
	and shrinking $I^\pm$ if necessary we shall assume that $\lambda'<B\left(t,\xi_0\right)$ for every $t \in \overline{I^+\cup I^-}$ (and we can no longer assume that $t_0\in\partial I^+$). Finally we rename $I^+$ and $I^-$ if necessary to conclude the following: there are $\lambda' \in \R$, a point $t_0\in S^1\setminus\left\{t_{\xi_0}\right\}$ and two open intervals $I^-\subset\left[t_{\xi_0},t_0\right], I^+\subset\left[t_0,t_{\xi_0}\right]$ such that
	\begin{enumerate}
		\item $\lambda'<B\left(t_0,\xi_0\right)$;
		\item $\sup_{t \in \overline{I\cup I'}} B\left(t,\xi_0\right) < \lambda'$.
	\end{enumerate}
	
	Since $m>0$, given $D>0$ there exists $n_0\in\N$ such that
	\[
	D\log\left|n\xi_0\right| < n^m\left(B\left(t_0,\xi_0\right)-\lambda'\right) \quad \mbox{and}\quad \left|I\right|,\left|I'\right|\ge \left|n\xi_0\right|^{-m}
	\]
	for every $n\ge n_0$. Then for each $t \in I^+\cup I^-$,
	\begin{align*}
		B\left(t_0,n\xi_0\right) - B\left(t,n\xi_0\right) &= n^m\left(B\left(t_0,\xi_0\right)-B\left(t,\xi_0\right)\right)\\
		&> n^m\left(B\left(t_0,\xi_0\right)-\lambda'\right)\\
		&> D\log\left|n\xi_0\right|, \quad \forall n\ge n_0.
	\end{align*}
	Therefore, for $n\ge n_0$ we conclude that $P$ does not satisfy $(\beta)_D$ at $n\xi_0$. Since $D$ is arbitrary, Theorem~\ref{Thm:main} shows that $P$ can not be globally solvable.
\end{proof}
\begin{Cor}
	\label{Cor:homog-charact}
	Let $c\left(t,D_x\right)$ be a positively homogeneous pseudodifferential operator of order $m\in\N$. Then $P=D_t+c\left(t,D_x\right)$ is globally solvable if and only if the following properties hold true:
	\begin{enumerate}
		\item $\left\{c_0\left(\xi\right)\right\}_{\xi\in\Z^N}$ satisfies~\eqref{eq:DC}.
		\item for each $\xi\in \Z^N\setminus\mathcal Z, \xi\ne0$, the function $b\left(t,\xi\right)$ does not change signal.
		\item for each $\xi\in\mathcal Z\setminus\{0\}$, given $\lambda \in \R$ the set
		\[
		\Omega_{\xi,\lambda} = \left\{t \in S^1 : B(t,\xi)<\lambda\right\}
		\]
		is connected.
	\end{enumerate}
\end{Cor}
\begin{proof}
	If $P$ is globally solvable then Theorem~\ref{Thm:main} ensures that $\left\{c_0\left(\xi\right)\right\}_{\xi\in\Z^N}$ satisfies~\eqref{eq:DC} while Proposition~\ref{Prop:homogeneous-complementar-Z-necessary} and Proposition~\ref{Prop:homogeneous-Z-necessary} ensure that the other conditions of the statement are satisfied. 
	
	Conversely, if $\xi \ne 0$ and $\xi \in \Z^N\setminus \mathcal Z$ then $c_0\left(\xi\right)\notin \Z$ and since $b(\cdot,\xi)$ does not change signal, $P$ satisfies $(\alpha)_1^+$ or $(\alpha)_1^-$ at $\xi$. Indeed, suppose that $b\left(t,\xi\right)\ge 0$ for every $t \in \left[0,2\pi\right]$. Then
	\begin{align*}
		B\left(t,\xi\right)-B\left(s,\xi\right) \le 0 < \log\left|\xi\right|, \quad \forall s \in \left[t,t+2\pi\right],
	\end{align*}
	so $P$ satisfies $(\alpha)_1^+$ at $\xi$. Analogously, one can show that $P$ satisfies $(\alpha)_1^-$ at any $\xi\in\Z\setminus\mathcal Z$, such that $\xi\ne 0$ and $b\left(t,\xi\right)\le 0$ for every $t \in \left[0,2\pi\right]$.

	Furthermore, given $\xi \in \mathcal Z\setminus\{0\}$, and $t \in S^1\setminus\left\{t_\xi\right\}$, the set
	\[
 	 \left\{s \in S^1 : B\left(s,\xi\right) \ge B\left(t,\xi\right)-\log\left|\xi\right|\right\} = S^1\setminus \Omega_{\xi,B\left(t,\xi\right)-\log \left|\xi\right|}
	\]
	is connected and contains $t$ and $t_\xi$, so it contains either $\left[t_\xi,t\right]$ or $\left[t,t_\xi\right]$, that is
	\[
	B\left(t,\xi\right) - B\left(s,\xi\right) \le \log \left|\xi\right|,\quad \forall t \in \left[t_\xi,t\right]
	\]
	or
	\[
	B\left(t,\xi\right) - B\left(s,\xi\right) \le \log \left|\xi\right|,\quad \forall t \in \left[t,t_\xi\right].
	\]
	This shows that $P$ satisfies $(\beta)_1$ at every $\xi\in\mathcal Z\setminus\{0\}$ and then $P$ is globally solvable by Theorem~\ref{Thm:main}.
\end{proof}

\section*{Acknowledgements}

I would like to express my gratitude to Professor Paulo D. Cordaro. Crucial points of this work would not have been achieved without our discussions and his valuable suggestions.

This work was supported by Fundação de Amparo à Pesquisa do Estado de São Paulo (FAPESP) under grant 2023/11769-5 and by Conselho Nacional de Desenvolvimento Cient{\'i}fico e Tecnol{\'o}gico (CNPq) under  grant 404175/2023-6. 
\bibliographystyle{plain}
\bibliography{bibliography}

\end{document}